\documentclass[reqno]{amsart}[12]
\usepackage{amssymb,amsmath,amsfonts,latexsym}
\usepackage{pb-diagram}
\usepackage[margin=1.25in]{geometry}
\numberwithin{equation}{section}

\newcommand{\complex}{{\mathbb{C}}}

\newcommand{\Z}{{\mathbb{Z}}}

\newcommand{\calC}{\mathcal{C}}
\newcommand{\calD}{\mathcal{D}}

\newcommand{\calO}{\mathcal{O}}

\newcommand{\Tr}{\operatorname{Tr}}
\newcommand{\tr}{\operatorname{tr}}
\newcommand{\Str}{\operatorname{Str}}

\newcommand{\End}{\operatorname{End}}

\newcommand{\gl}{\mathfrak{gl}}

\renewcommand{\det}{\operatorname{det}}

\theoremstyle{plain}
        \newtheorem{theorem}{Theorem}[section]
        
        \newtheorem{proposition}[theorem]{Proposition}
        \newtheorem{corollary}[theorem]{Corollary}
        
\theoremstyle{definition}
        
        \newtheorem{remark}[theorem]{Remark}

\title{Equivariant Lefschetz number of differential operators}
\date{\today}
\author{G. Felder \textrm{and} X.~Tang}
\begin{document}
\maketitle
\begin{abstract}
Let $G$ be a compact Lie group acting on a compact complex
manifold $M$. We prove a trace density formula for the
$G$-Lefschetz number of a differential operator on $M$. We
generalize  Engeli and Felder's recent results to orbifolds.
\end{abstract}

\section{Introduction}
In this paper, we study a $G$-equivariant Lefschetz
number formula for a compact Lie group action on a compact
connected complex manifold $M$. We assume that the isotopy group at each  point of $M$ is finite, which implies that the quotient space $M/G$ is an orbifold.

Let $E$ be a $G$-equivariant vector bundle on $M$. Both $G$ and a
global holomorphic differential operator $D\in \calD_E(M)$ act on
the sheaf cohomology group $H^i(M,E)$. Since $M$  is compact, the
sheaf cohomology $H^i(M,E)$ is finite dimensional. We can consider
the following $G$-Lefschetz number, for $\gamma\in G$
\[
D\mapsto L(\gamma,D)=\sum_i(-1)^i\tr(\gamma H^i(D)),
\]
where $H^i(D)$ denotes the action of $D$ on $H^i(M, E)$.
$L(\cdots, D)$ is a smooth function on $G$ which
contains useful information of the $G$-action. When $D$ is the
identity operator,  the topological
expression of this number is known as the $G$-equivariant Riemann-Roch-Hirzebruch  formula \cite{BeGeVe} for a $G$-equivariant bundle $E$.\\

In this paper, we express the number $L(\gamma, D)$ for any
differential operator, by an integral over the $\gamma$ fixed
point submanifold $M^\gamma$  of some differential form
$\chi_{0,\gamma}$ determined by finitely many jets of $D$ on
$M^g$ and a hermitian metric on $M$ and $E$,
\[
L(\gamma,D)=\int_{M^\gamma}\frac{1}{(2\pi
i)^{n-l(\gamma)/2}}\chi_{0,\gamma}(D),
\]
where $l(\gamma)$ is the real codimension of $M^\gamma$ in $M$.
The proof of this theorem is a generalization of Engeli-Felder's
theorem \cite{enfe} on a manifold, with the local results
developed in the previous work \cite{ppt}. The analysis of this
proof is more involved than \cite{enfe} because we have to work
with the $\gamma$ traces of a heat kernel. When $G$ is finite and
$D$ is $G$ invariant, we are  able to compute the $G$ average of
$L(\gamma,D)$,
\begin{equation}\label{eq:average}
\frac{1}{|G|}\sum_{\gamma \in G} L(\gamma,
D)=\sum_{(\gamma)\subset
G}\int_{M^\gamma/C(\gamma)}\frac{1}{m_\gamma}\frac{1}{(2\pi
i)^{n-l(\gamma)/2}}\chi_{0,\gamma}(D),
\end{equation}
where  $(\gamma)$ stands for the conjugacy class in $G$ containing
$\gamma$, $C(\gamma)$ is the centralizer group of $\gamma$ in $G$,
and $m_\gamma$ is the number of isotopy of $C(\gamma)$ action on
$M^\gamma$. The proof of this result is an application of the
proof of the $G$-Lefschetz number
formula.\\

A $G$ invariant differential operator $D$ on $M$ descends to a differential operator $\hat{D}$ on the orbifold $X=M/G$. Hence, the above equation (\ref{eq:average}) can be viewed as a Lefschetz number formula of $\hat{D}$.  This inspires a question of computing the Lefschetz number of a general differential operator on an orbifold $X$. It turns out that there are at least two in-equivalent definitions of a differential operator on an orbifold. One definition of a differential operator on $X$ from a vector bundle $E$ to $F$ is a linear map of sections $\Gamma^\infty(E)\to \Gamma^\infty(F)$ expressed by a finite combination of bundle endomorphisms and their covariant derivatives; the other definition is that a linear map $P:\Gamma^\infty(E)\to \Gamma^\infty(F)$ such that there is an integer $N$ and for any smooth functions $f_0, \cdots, f_N$ on $X$, (holomorphic functions if one considers a holomorphic differential operator), the commutator
\[
[f_N,[f_{N-1},\cdots, [f_0, P]\cdots]]=0.
\]
In the case of a manifold, these two definitions are equivalent.
However, in the case of an orbifold, it is quite easy to check
that the second definition contains the first one, but not vice
versa. We provide in Remark \ref{rmk:diff} an explicit example
which is a differential operator in the sense of the second
definition but not the first on the orbifold
$\complex/\mathbb{Z}_2$. In the following, we call  operators in
the first definition geometric differential operators, and those
in the second definition algebraic ones.

We establish a Lefschetz number formula for a geometric differential operator on  a general (maybe non-reduced) complex orbifold. Let $E$ be
an vector bundle on a compact orbifold $X$. We consider a global geometric differential operator $D$ on $X$, which is a global section of the sheaf of geometric
holomorphic differential operators acting on sections $E$.  Since
$X$ is compact, we have that $H^i(X,E)$ is finite dimensional, and the following orbifold Lefschetz number is well defined,
\[
D\longmapsto L(D)=\sum_{i}(-1)^i\tr(H^i(D)).
\]
We prove in this paper an integral formula to compute $L(D)$, i.e.
\[
L(D)=\int_{\widetilde{X}}\frac{1}{m_\calO}\frac{1}{(2\pi i)^{n-l(\calO)/2}}\chi_{0, \calO}(D),
\]
where $\widetilde{X}$ is the inertia orbifold associated to the
orbifold $X$, and $l(\calO)$ is a local constant, the real
codimension of $\widetilde{X}$ in $X$, and $m_\calO$ is also a
local constant, the number of isotopy, and $\chi_{0,\calO}(D)$ is
a top differential form on $\widetilde{X}$.   The proof of this
result is an extension of the proof of the $G$-equivariant
Lefschetz number formula. When $D$ is the identity operator, then
the above theorem together with the local Riemann-Roch-Hirzebruch
theorem in \cite{ppt} computes the holomorphic Euler
characteristic of $E$ on $X$, which was computed in
\cite{du:book}. \\

We remark that the $G$ Lefschetz formula proved in this paper
probably can be generalized and applied to study toric varieties.
For a compact toric manifold, the canonical torus action does not
always have finite isotopy groups but only has finitely many
isolated fixed points and each fixed point comes with a nice
coordinate system. Musson \cite{mu:toric} provided a description
of the algebra of differential operators on toric varieties. We
plan to study the Lefschetz number formula and its application on
toric varieties in the future.

This paper is designed as follows. In Section 2, we briefly review
the $\gamma$-twisted Hochschild (co)homology of differential
operators and some formal differential geometry on an orbifold. In
Section 3, we prove the Lefschetz number formula for $L(\gamma,
D)$ for $G$ acting on a compact manifold $M$. In Section 4, we
compute the $G$-average and orbifold Lefschetz number formula. And
we end this paper with a remark on differential operators on
orbifolds and an open
question for future research.\\

\noindent{\bf Acknowledgments:} We would like to thank the
organizers of the trimester on ``Groupoids and stacks in physics
and geometry" to hosting their visits of Institut Henri
Poincar\'e. This work has been partially supported by the MISGAM
programme of the European Science Foundation. The first author was
partially supported by the European Union through the FP6 Marie
Curie RTN ENIGMA (Contract number MRTN-CT-2004-5652) and the Swiss
National Science Foundation (grant 200020-105450). The second
author was partially supported by the U.S. National Science
Foundation (grant 0703775).

\section{Hochschild cohomology and formal geometry}
In this section, we review some known results about Hochschild
(co)homology of the algebra of differential operators on an
orbifold. In this paper, we assume that $G$ is a compact group acting on a complex manifold $M$ such that the isotopy group at each point of $M$ is finite.

Let $\gamma$ be an element of $G$. We study the geometry near a
fixed point $x$ of $\gamma$. By Bochner's  theorem \cite{bochner},
nearby $x$, there is a $\gamma$ invariant coordinate neighborhood
$U$ such that the $\gamma$ action is linear.  And because of the
above  assumption of finite isotopy, $\gamma$ acts on $U$ of
finite order. This inspires the following consideration.

Let $\gamma$ be an element of a compact group acting on an
$n$-$\dim$ complex vector space. We consider
$\calO_n=\complex[[y_1, \cdots, y_n]]$ and
$\calD=\calO_n[\partial_{y_1}, \cdots,
\partial _{y_n}]$ the algebra of formal
power series in $n$ variables and formal differential operators,
and $\calO^{\text{pol}}_n=\complex[y_1, \cdots, y_n]$ and
$\calD^{\text{poly}}=\calO^{\text{pol}}_n [\partial_{y_1}, \cdots,
\partial_{y_n}]$ their subalgebras of polynomial functions
and differential operators. The $\gamma$-twisted Hochschild
homology $H_\bullet(\calD^{\text{pol}},
\calD_\gamma^{\text{pol}})$ of $\calD^{\text{pol}}$ is computed to
be (\cite{AFLS:HIAWAGF})
\[
\left\{\begin{array}{ll}\complex&\bullet=2n-l(\gamma),\\ 0&\text{else}. \end{array}\right.
\]
where $l(\gamma)$ is the real codimension of the $\gamma$-fixed
point subspace $V^\gamma$. And by $\calD^{\text{pol}}_\gamma$, we
mean the vector space $\calD^{\text{pol}}$ with the following
bimodule structure
\[
f\cdot\xi=f\circ\xi,\ \ \ \ \ \ \ \xi\cdot
f=\xi\circ\gamma^{-1}(f).
\]

We decompose $\complex^n$ into
$V^\gamma\oplus V^\perp$ where $V^\perp$ is the $\gamma$ invariant subspace of
$\complex^n$ complement to $V^\gamma$. Let  $y^{\gamma}_1,\cdots, y^{\gamma}_{l(\gamma)/2}$ be a basis of $V^\gamma$. Then the following expression defines a generator for
$H_{2n-l(\gamma)}(\calD, \calD_\gamma)$,
\[
c^{\gamma}_{2n-l(\gamma)}=\sum_{\epsilon\in S_{2n-l(\gamma)}}1\otimes
u_{\epsilon(1)}\otimes\cdots u_{\epsilon(2n-l(\gamma))},\ \ \
u_{2i-1}=\partial_{y^{\gamma}_i},\ u_{2i}=y_{i}^\gamma.
\]

An explicit Hochschild cocycle $\tau_{2n-l(\gamma)}^{\gamma}$ on
$H^\bullet(\calD; \calD^*_{\gamma})$ was constructed in \cite{ppt}
as follows. Given any $D\in \calD$, we can decompose it in $\sum_i
D^\gamma _i\otimes D^\perp_{i}$ with $D_i$ a differential operator
on $V^\gamma$ and $D^\perp_i$ a differential operator on
$V^\perp$. Define
\[
\begin{split}
&\tau_{2n-l(\gamma)}^{\gamma}(D_0,
\cdots, D_{2n-l(\gamma)})\\
=&\sum_{i_0, \cdots,
i_{2n-l(\gamma)}}\tau_{2n-l(\gamma)}(D^\gamma_{0, i_0}, \cdots,
D^\gamma_{2n-l(\gamma), i_{2n-l(\gamma)}})\Tr_\gamma(D^\perp_{0,
i_0}\circ\cdots \circ D^\perp_{2n-l(\gamma), i_{2n-l(\gamma)}}),
\end{split}
\]
where we have written $D_k=\sum_{i_k}D^\gamma_{k, i_k}\otimes
D^\perp_{k, i_k}$, and $\tau_{2n-l(\gamma)}$ is a cocycle in
$H^{2n-l(\gamma)}(\calD_{V^\gamma}, \calD_{V^\gamma}^*)$ as is
defined in \cite{FFS}, and $\Tr_\gamma$ is the $\gamma$-trace on
$\calD_{V^\perp}$ defined in \cite{fe:g-index}. We refer to
\cite{ppt}[Section 3] for explicit formulas of the cocycles
$\tau_{2n-l(\gamma)}$ and $\Tr_\gamma$.\\

The above computation of Hochschild (co)homologies can be extended
to matrix valued differential operators $\calD_E$ with $E$ a
finite dimensional representation of $\gamma$. The cycle
$c^\gamma_{E, 2n-l(\gamma)}$ and cocycle $\tau^{\gamma}_{E,
2n-l(\gamma)}$ can be defined to be  cocycles on $\calD_E$.

We list some properties of $\tau^{\gamma}_{E, 2n-l(\gamma)}$,
which were proved in \cite{ppt}[Section 3].
\begin{enumerate}
\item $\tau^{\gamma}_{E, 2n-l(\gamma)}$ is invariant under the
action of $(GL_n(\complex)\times GL(E))^\gamma$, the subgroup of
$GL_n(\complex)\times GL(E)$ consisting elements commuting with
$\gamma$;
\item $\tau^{\gamma}_{E, 2n-l(\gamma)}(c^{\gamma}_{E, 2n-l(\gamma)})=(\det(1-\gamma^{-1}))^{-1}\tr_E(\gamma)$, where we
consider the determinant of $1-\gamma^{-1}$ on $V^\perp$.
\end{enumerate}

\begin{remark}
From the above formula it is possible that if $\tr_E(\gamma)=0$, $\tau_{E, 2n-l(\gamma)}^\gamma(c_{E,2n-l(\gamma)})=0$. However, $c_{E, 2n-l(\gamma)}^\gamma$ of the same form as $c^\gamma_{2n-l(\gamma)}$ is never 0, which makes the later arguments work.
\end{remark}

The following theorem is a straight forward generalization of
Brylinski-Getzler's results to the $\gamma$-twisted case, whose
proof is omitted.

\begin{theorem}
\label{thm:local}Let $\gamma$ be an element of a compact group
acting on $M$ preserving its complex structure. At every $\gamma$
fixed point $x$, there is a $\gamma$ invariant coordinate
neighborhood $U$ on which $\gamma$ acts linearly.  Let $E$ be a
$\gamma$-equivariant vector bundle on $M$. On $U$,
$H_\bullet(\calD_E(U), {\calD_E(U)}_\gamma)$ is spanned by $c_{E,
2n-l(\gamma)}^{\gamma}$.
\end{theorem}

We consider the setting in Theorem \ref{thm:local} and a connected
component of $\gamma$ fixed point, which is denoted by $M^\gamma$.
Let $E$ be a $\gamma$-equivariant complex vector bundle on $M$,
which induces a bundle $E$ on $M^\gamma$. On $M^\gamma$, there is
a vector bundle $N^\gamma$ which is the normal bundle associated
to the embedding of $M^\gamma$ in $M$. We notice that $\gamma$
acts on both $N^\gamma$ and $E$ fiberwisely. Let $\calO_{N,\gamma}(E)$ be a bundle over $M^\gamma$ whose fiber at $x\in M^\gamma$ is the algebra of $\End(E_x)$ valued
holomorphic functions on $N_x^\gamma$. We consider the
sheaf of Lie algebra of $\gamma$-invariant vector fields, i.e.
\[
W_n^\gamma=(\oplus_i\calO_n\partial_{y_i})^\gamma,\ \text{and}\
W_n(E_x)^\gamma=\big( (\oplus_i\calO_n\partial_{y_i})\ltimes
\gl_r(\calO_n(E_x))\big)^\gamma,\ \ r=rank(E).
\]
The formal infinite jets $J^\infty(\calO_{N,\gamma}(E))$ of $\calO_{N,\gamma}(E)$ is
a principal $W_n(E)^\gamma$-space. Furthermore, there is a natural
flat connection $A_{E,\gamma}$ on $J^\infty(\calO_{N,\gamma}(E))$ with
\[
dA_{E,\gamma}+\frac{1}{2}[A_{E,\gamma}, A_{E,\gamma}]=0.
\]

The connection $A_{E,\gamma}$ induces a flat connection on the
associated bundle $J_1(E)\times_G \calD_{\calO_{N,\gamma}(E)}$
over $M^\gamma$, where $G=(GL_{n-l(\gamma)/2}\times GL_r)^\gamma$,
with $r=rank(E)$. The flat sections of this bundle $J_1(E)\times_G
\calD_{\calO_{N,\gamma}(E)}$ are in one-to-one correspondence to
the restriction to $M^\gamma$ of global differential operators
$\calD_E$ on $M$. We denote this flat connection on
$J_1(E)\times_G \calD_{\calO_{N,\gamma}(E)}$ by
$\omega_{E,\gamma}$.

\begin{proposition}\label{prop:local-hoch}
Let $\Omega^\bullet$ be the $l(\gamma)$-shifted complex of sheaves of
complex valued smooth differential forms on $M^\gamma$ with
the de Rham differential, and $\calC_{\bullet, \gamma}(\calD)$ be the
$\gamma$-twisted complex of Hochschild chains of $\calD$. There is a
homomorphism of complexes of sheaves on $M^\gamma$
\[
\chi_\gamma: \calC_{\bullet,\gamma}(\calD)\rightarrow
\Omega^{2n-\bullet}(M^\gamma)[l(\gamma)].
\]
\end{proposition}

\begin{proof}
The proof of this proposition is a copy of the proof of
Proposition 2.4.1 \cite{enfe}. We have the following sequence of
maps
\[
\chi_{p,\gamma}: \calC_{p,\gamma}(\calD)\rightarrow \Omega^{2n-l(\gamma)-p}(M^\gamma)
\]
by
\[
\chi_{p,\gamma}(D_0, \cdots, D_p)=(-1)^p\tau^\gamma_{V, 2n-l(\gamma)}(\text{sh}_{p,2n-l(\gamma)-p}
(\hat{D}_0, \hat{D}_1,\cdots, \hat{D}_p, \omega_{E,\gamma}, \cdots, \omega_{E,\gamma})),
\]
where $\text{sh}_{p,2n-l(\gamma)-p}$ is the sum over all $(p,
2n-l(\gamma)-p)$ shuffles, and we have identified a differential
operator $D_i$ with a flat section $\hat{D}_i$ of the previously
considered bundle $J_1(E)\times_G \calD_{\calO_{N,\gamma}(E)}$.

It is straightforward to check the following identities
\[
d\circ \chi_{p,\gamma}=\chi_{p-1,\gamma}\circ b^\gamma.
\]
And the proposition follows from this identity.\\
\end{proof}

We conclude this section by introducing two $\gamma$-traces on $\calD_E$. In the following section, we will prove that these two traces are actually equal.

\begin{enumerate}
\item The first trace $\Tr_{1,\gamma}:\calD_E(M)\to \complex$ is defined to be
\[
D\longmapsto \sum_{j=0}^n(-1)^j\tr_j(\gamma H^j(D)),
\]
where $H^j(M,E)$ is the sheaf cohomology group with the induced
$\gamma$ action, and $H^j(D):\calD_E(M)\to \End(H^j(M, E))$. We
remark that as $M$ is compact, $H^j(M,E)$ is finite dimensional,
and the number $\Tr_{1,\gamma}$ is a finite number.
\item The second trace $\Tr_{2, \gamma}$ is defined to be
\[
D\longmapsto \int_{M^\gamma}\frac{1}{(2\pi
i)^{n-\frac{l(\gamma)}{2}}}\chi_{0,\gamma}(D).
\]
We remark that in the above formula, if $M^\gamma$ has different components, we sum over all components.  And because $M$ is compact,  each component of $M^\gamma$ is compact and there are only finitely many components, and therefore the above integral is a finite number.
\end{enumerate}
\section{G-Lefschetz number}
We prove in this section that $\Tr_{1, \gamma}=\Tr_{2, \gamma}$ by
verifying that they are both equal to a third trace which will be
introduced in this section. Many ideas of the proof in this
section are originally from \cite{enfe}. We adapt them to study
the $\gamma$-twisted traces.

In this whole section, we assume that $\gamma$ is an element of a
compact group $G$ acting on $M$ preserving the complex structure
and the isotopy of $G$ action at each point of $M$  is finite.

\subsection{The third trace}
We introduce a third $\gamma$-twisted trace on $\calD_E(M)$ in
this subsection, and will prove in the following subsections that
it is equal to the first and second $\gamma$-twisted traces
separately.\\

We choose $(U_i)$ a locally finite open cover of $M^\gamma$ and
consider the sheaf $\calD_{\calO_{N,\gamma}(E)}$ of differential
operators on the bundle $\calO_{N,\gamma}(E)$. We remark that flat
global sections of this sheaf are one-to-one correspondent to
restrictions of global differential operators $\calD_E$ to
$M^\gamma$. We consider the \v{C}ech double complex
$\check{C}^q(U, \calC_{-p}(\calD_E,\calD_{E,\gamma}))$, where
$\calC_\bullet(\calD_{E},\calD_{E,\gamma})$ is the complex of
$\gamma$ twisted Hochschild chains on
$\calD_{\calO_{N,\gamma}(E)}$.

Given any operator $D\in \calD_E(M)$, the restriction of $D$ to
$M^\gamma$ defines a $(0,0)$ cocycle in $\check{C}^0(U, \calC_0)$.
By Theorem \ref{thm:local}, $D|_{U_i}$ is a Hochschild boundary
when restricted to a sufficiently small open set $U_i$. Therefore,
there is a $(-1,0)$-cochain $D^{(1)}$ in $C^{-1,0}$ such that
$D|_{U_i}=b^\gamma D^{(1)}|_{U_i}$. Thereafter, we consider the
\v{C}ech differential on $D^{(1)}$, i.e.
$\delta(D^{(1)})|_{ij}=D^{(1)}|_{U_i}-D^{(1)}|_{U_j}$. As $\delta$
commutes with $b^\gamma$,
$b^\gamma(\delta(D^{(1)}))=\delta(b^{\gamma}(D^{(1)}))=\delta(D)=0$.
This shows that $\delta(D^{(1)})$ is again a Hochschild cycle. And
again by Theorem \ref{thm:local}, we know that there is an element
$D^{(2)}$ in $C^{-2,1}$ such that
$b^{\gamma}D^{(2)}=\delta(D^{(1)})$. By repeating this induction
step, we will have a sequence of cochains $D^{(j)}\in C^{-j,
j-1},\ j=1,\cdots, 2n-l(\gamma)$ with
\[
b^\gamma D^{(1)}=D, \ \ \ \delta(D^{(j)})=b^\gamma(D^{j+1}),\ \
j=1,\cdots, 2n-l(\gamma)-1.
\]

We need to keep in mind at the step $j=2n-l(\gamma)$ the
Hochschild homology is not trivial. Therefore, we have that
\[
\delta(D^{(2n-l(\gamma))})=s^\gamma+b^\gamma(D^{(2n+1-l(\gamma))}),
\]
where $s^\gamma\in C^{2n-l(\gamma), -2n+l(\gamma)}$ is equal to
\begin{equation}\label{eq:s}
s^\gamma_{i_0, \cdots, i_{2n-l(\gamma)}}=\lambda^\gamma_{i_0,
\cdots, i_{2n-l(\gamma)}}(D)c^\gamma_E(U_{i_0}\cap\cdots\cap
U_{i_{2n-l(\gamma)}}),
\end{equation}
for some $2n-l(\gamma)$ \v{C}ech cocycle $\lambda^\gamma(D)$ in
$\check{C}$$^{2n-l(\gamma)}(U;\complex)$. Therefore,
$\lambda^\gamma(D)$ is an element in $H^{2n-l(\gamma)}(M^\gamma)$.
We define $\Tr_{3,\gamma}(D)$ to be the Poincar\'e dual of
$[\lambda^\gamma(D)]$, which is a number in $\complex$.\\

\subsection{Local expression of the trace density}
We start by fixing a smooth and locally finite triangulation
$|K^\gamma|$ to $M^\gamma$, with the underlying simplicial complex
$K_0^\gamma$. We consider an open cover $(U_i)_{i\in K_0^\gamma}$
of $M^\gamma$, where $K_0^\gamma$ is the set of vertices of the
triangulation, such that $U_i$ is the complement of the simplexes
not containing the vertex $i$. Because of the construction, we see
that the cover $(U_i)$ satisfies for any $i_1<\cdots <i_p$,
\begin{enumerate}
\item $U_{i_0}\cap \cdots \cap U_{i_p}$ is contractible,
\item If $p\gg2n-l(\gamma)$, then $U_{i_0}\cap\cdots \cap U_{i_p}$
is empty.
\end{enumerate}

We consider the cell decomposition $C_\bullet$ dual to the above
triangulation $K^\gamma$. We denote $C_{i_0,\cdots, i_p}$ the
$(2n-l(\gamma)-p)$-cells dual to the simplex $K_{i_0, \cdots,
i_p}$. The orientation of $C_\bullet$ is set to require $C_{i_0,
\cdots, i_p}\cdot K_{i_0, \cdots, i_{p}}=1$.

We have the following analogous proposition as
\cite{enfe}[Prop. 5.1].
\begin{proposition}\label{prop:trace-2}
For any global differential operator $D$ on $M$. Let $s$ be the
$(2n-l(\gamma), -2n+l(\gamma))$-$\check{C}$ech cocycle defined in
Equation (\ref{eq:s}), and $\chi_{2n-l(\gamma), \gamma}$ be the map
defined in Proposition \ref{prop:local-hoch}. Then
\[
\Tr_{2,
\gamma}(D)=(-1)^{n-l(\gamma)/2}\sum_{i_0<\cdots<i_{2n-l(\gamma)}}\int_{C_{i_0,\cdots,
i_{2n-l(\gamma)}}}\chi_{2n-l(\gamma), \gamma}(s^\gamma_{i_0,
\cdots, i_{2n-l(\gamma)}}).
\]
\end{proposition}
\begin{proof}
The proof is same to the proof of Proposition 5.1, \cite{enfe}.
\end{proof}
With this proposition, we can compare $\Tr_{2, \gamma}$ and
$\Tr_{3,\gamma}$ locally on each $C_{i_1, \cdots,
i_{2n-l(\gamma)}}$ and verify directly that on each connected component
\[
\Tr_{2,\gamma}=(-1)^{n-l(\gamma)}\frac{\tr_E(\gamma)}{\det(1-\gamma^{-1})}\sum_{i_0<\cdots
<i_{2n-l(\gamma)}}\lambda^\gamma_{i_0,\cdots,
i_{2n-l(\gamma)}}(D)=(-1)^{n-l(\gamma)}\frac{\tr_E(\gamma)}{\det(1-\gamma^{-1})}\Tr_{3,\gamma}(D).
\]

\begin{remark}
Because $\gamma$ acts on the bundle $E|_{M^\gamma}$ and $N^\gamma|_{M^\gamma}$ fiberwisely of finite order, the eigenvalues of $\gamma$ are discretely distributed. This implies that $\tr_E(\gamma)$ and $\frac{1}{1-\gamma^{-1}}$ are both local constants on $M^\gamma$, and the above equation is well defined.
\end{remark}
\subsection{Asymptotic pairing}
This subsection provides an important tool to prove in the
next subsection that $\Tr_{1,\gamma}=\Tr_{3,\gamma}$.\\

We consider the Dolbeault complex $(\Omega^{(0,\bullet)}(M,E),
\bar{\partial})$ with values in the holomorphic vector bundle $E$.
Because of the properness assumption on $\gamma$ action, we can fix a
$\gamma$-invariant hermitian metric on $T_M$ and $E$. Accordingly,
we can consider the Hilbert space of $L^2$-integrable sections of
$\Omega^{(0,\bullet)}(M,E)$. On this Hilbert space, there is a
self-adjoint positive semidefinite operator
$\Delta_{\bar{\partial}}=\bar{\partial}\bar{\partial}^*+\bar{\partial}^*\bar{\partial}$,
where $\bar{\partial}^*$ is the Hodge dual of $\bar{\partial}$.

Let $e^{-t\Delta_{\bar{\partial}}}$ be the heat operator with
kernel $k_t$. According to \cite{BeGeVe}[Theorem 2.30], as $t\to 0$ the heat
kernel has an asymptotic expansion,
\begin{equation}\label{eq:expansion}
k_t(z,z')\sim\frac{1}{(\pi
t)^n}e^{-\frac{dist(z,z')^2}{t}}(\Phi_0(z,z')+t\Phi_1(z,z')+\cdots),
\end{equation}
where $dist(z,z')$ is the geodesic distance between $z$ and $z'$.

We have the following generalization of Proposition 6.1 in
\cite{enfe}.
\begin{proposition}\label{prop:jlo}Let $U$ be a $\Gamma$-invariant open subset of
$M$, $A=\calD_E(U)$, and
$M_c=\Omega^{0,\bullet}_c(U)\otimes_{\calO_M(U)}\calD_E(U)$. We
consider the $k_t^N$ be the truncation at the $N$-th term of the
asymptotic expansion (\ref{eq:expansion}) with support in a small
neighborhood of the diagonal in $U\times U$. For $D_0\in
M_c^k=\Omega^{0,k}_c(U)\otimes_{\calO_{M}(U)}\calD_E(U)$,
$D_1,\cdots, D_k\in A$, the expression
\begin{equation}\label{eq:jlo}
\Psi^\gamma_k(D_0, \cdots,
D_k)=(-1)^{\frac{k(k+1)}{2}}\big[\int_{t\Delta_k}\Str(\gamma
D_0k_{s_0}^N[\bar{\partial}^*, D_1]k_{s_1}^{N}\cdots
[\bar{\partial}^*, D_k]k_{s_k}^N)ds_1\cdots ds_k\big]_-
\end{equation}
is independent of $N$ for $N\gg 1$ and defines a continuous
cocycle
\[
\Psi^\gamma=\sum_k\Psi^\gamma_k\in \bigoplus_{k=0}^n
Hom(M_c^k\otimes\bar{A}^{\otimes k},
\complex)[t^{-\frac{1}{2}}]=C^0(A,
M_{c,\gamma}^*)[t^{-\frac{1}{2}}].
\]
In the above formula for $\Psi^\gamma_k$, $\Str$ is the super trace
on $\wedge^\bullet {T^{0,1}}^*U\otimes E|_U$, and $[\cdots ]_-$
takes the none positive $t$-power terms.
\end{proposition}

\begin{proof}
Once we notice that the appearance of $\gamma$ in $\Psi^\gamma$
leads to the twisted cocycle condition and the fact that
$\bar{\partial}, \bar{\partial}^*, \Delta_{\bar{\partial}}$ are
all $\gamma$-invariant, the proof of this proposition is a repeat
of the proof of Proposition 6.1 in \cite{enfe}.
\end{proof}

In the following, we adapt Engeli-Felder's construction in \cite{enfe}[Section 6] to the
$\gamma$-twisted situation.

Given $\varphi_0,\cdots, \varphi_k\in C_c^\infty(U)\subset A$, we
can view them as $0$-cochains in $C^\bullet(A, M_{c})$. Let $\delta$ be the
differential on $C^\bullet(A, M_{c})$. We consider
\[
Z^k=\varphi_0\cup \delta\varphi_1\cup \cdots \delta\varphi_k\in C^k(A, M_{c}),
\]
where $\cup:C^p(A;M_{c})\otimes
C^q(A;M_{c})\rightarrow C^{p+q}(A;M_{c}\otimes_A
M_{c})= C^{p+q}(A;M_{c})$ by $M_{c}\otimes_A
M_{c}\to M_{c}$.

By
$M_{c,\gamma}$, we mean the linear space $M_c$, but the right
action of $A$ is twisted by $\gamma$. We use the following cup product
\[
\cup:C^\bullet(A; M_{c,\gamma}^*)\otimes
C^\bullet(A;M_c)\rightarrow C^{\bullet}(A; M_{c,\gamma}^*\otimes_A
M_c)=C^\bullet(A; A^*_\gamma),
\]
to construct an element $\sigma^\gamma_k(\varphi_0, \cdots, \varphi_k)$ in $C^\bullet(A;A^*_\gamma)$ by
\[
\sigma^\gamma_k(\varphi_0, \cdots, \varphi_k)=\Psi^\gamma\cup Z^k(\varphi_0, \cdots, \varphi_k).
\]

If $A=\calD_E(U)$, $B=C^\infty(U)$, the above $\sigma^\gamma_k$
defines a linear map
\[
\sigma^\gamma_k:C_k(A;A_\gamma)\otimes C_k^c(B)\rightarrow \complex[t^{-\frac{1}{2}}].
\]

We have the following properties of the map $\sigma^\gamma_k$.
\begin{proposition}
\label{prop:pairing}
\noindent{}
\begin{enumerate}
\item $\sigma_k^\gamma$ vanishes on $(\varphi_0, \cdots, \varphi_k)$ with
$\cap_{i=0}^k \text{supp}( \varphi_i)=\varnothing$;

\item If $\varphi=\varphi_0\otimes \cdots \otimes
\varphi_k$, $s(\varphi)=1\otimes \varphi_0\otimes\cdots\otimes
\varphi_k$,  and $D\in C_{k+1}(A;A_\gamma)$, then
$\sigma^\gamma_k(b^\gamma D\otimes \varphi)=
\sigma^\gamma_{k+1}(D\otimes s(\varphi))$ for $k\geq 0$;

\item $\sigma_0^\gamma(D,\varphi)=\big[\sum_{j=0}^{2n-l(\gamma)}\tr_{\Omega^{0,j}}
(\varphi\gamma De^{-t\Delta_{\bar{\partial}}})\big]_-$;

\item If $\varphi_i$,
$i=1, \cdots, 2n-l(\gamma)$ are $\gamma$ invariant and constant
along the normal directions of $U^\gamma$ within a tubular
neighborhood of $U^\gamma$,
$\sigma^\gamma_{2n-l(\gamma)}(c^\gamma_E(U);
\varphi_0\otimes\cdots\otimes \varphi_{2n-l(\gamma)})$ is equal to
\[
\frac{(-1)^{n-l(\gamma)/2}\tr_E(\gamma)}{(2\pi i)^{n-l(\gamma)/2}}
\int_{U^\gamma}\frac{1}{\det(1-\gamma^{-1})}\varphi_0d\varphi_1\cdots
d\varphi_{2n-l(\gamma)}.
\]
\end{enumerate}
\end{proposition}
\begin{proof}
We write the pairing $\sigma_k^\gamma$ in a more explicit way:
\begin{equation}\label{eq:pairing}
\sigma^\gamma_k(D_0, \cdots, D_k; \varphi_0, \cdots,
\varphi_k)=\sum_{j=0}^k
(-1)^{j(k-j)}\Psi^\gamma_j(\gamma^{-1}(Z^k_{k-j}(D_{j+1}, \cdots,
D_k; \varphi_0, \cdots, \varphi_k))D_0, \cdots, D_j),
\end{equation}
where $Z^k_{k-j}(D_{j+1}, \cdots, D_k; \varphi_0, \cdots, \varphi_k)$ is equal to
\[
\sum_{\pi\in S_{k-j,j}}\text{sign}(\pi)\varphi_0B_{\pi(1)}(\varphi_1)\cdots B_{\pi(k)}(\varphi_k),
\]
with $B_i(\varphi)=[D_{j+i}, \varphi]$ for $i=1, \cdots, k-j$,
and $B_i(\varphi)=[\bar{\partial}, \varphi]$, for $i=k-j+1, \cdots, k$.\\

For (1), we see that if
$\cap_{i=0}^k\text{supp}(\varphi_i)=\varnothing$, then
$Z^k_j(\varphi_0, \cdots, \varphi_k)=0$ everywhere. This implies
$\sigma_k^\gamma$ vanishes on $(\varphi_0,\cdots, \varphi_k)$.\\

For (2), we notice that $\delta Z^k(\varphi_0, \cdots,
\varphi_k)=\delta\varphi_0\cup \cdots \cup \delta\varphi_k=1\cup
\delta \varphi_0\cup \cdots \cup \delta \varphi_k=Z^{k+1}(1,
\varphi_0, \cdots, \varphi_k)$. Hence, since $\Psi^\gamma$ is a
$\gamma$-twisted cocycle, we have $b^\gamma(\Psi^\gamma\cup
Z(\varphi))=\Psi^\gamma\cup \delta(Z(\varphi))=\Psi^\gamma\cup
Z(s(\varphi))$.\\

For (3), we can check it directly by definition using the above
explicit formula (\ref{eq:pairing}) for $\sigma^\gamma$.\\

For (4), we recall that $c^\gamma_{2n-l(\gamma)}=1\otimes u$,
where $u$ contains $(2n-l(\gamma))!$ terms with $\partial_{y^i}$
or $y^i$. We observe from Equation (\ref{eq:pairing}) and the
definition of $Z^{2n-l(\gamma)}_j$ that if the multiplication
operator by $y^i$ appears in $Z^{2n-l(\gamma)}_j$, then
$Z_j^{2n-l(\gamma)}$ vanishes because the commutator $[y^i,
\varphi]=0$. This implies that the none zero terms in
$\sigma^\gamma_{2n-l(\gamma)}$ contains only terms in the
expression (\ref{eq:pairing}) with $j\geq n-\frac{l(\gamma)}{2}$.

On the other hand, as $\varphi_i$ is independent of the direction
along the normal direction to $U^\gamma$ within a tubular
neighborhood of $U^\gamma$, $[\bar{\partial}, \varphi_i]$ is an
anti-holomorphic differential form along the direction of
$U^\gamma$. As $\dim_\complex(U^\gamma)=n-l(\gamma)/2$, inside
this tubular neighborhood of $U^\gamma$,
$Z_{2n-l(\gamma)-j}^{2n-l(\gamma)}(\cdots)$ cannot contain more
than $n-l(\gamma)/2$ terms like $[\bar{\partial}, \varphi_i]$,
because otherwise $Z_{2n-l(\gamma)-j}^{2n-l(\gamma)}(\cdots)$
contains a product of more than $n-l(\gamma)/2$ terms of
anti-holomorphic differential forms along $U^\gamma$. This implies
that when $j>n-\frac{l(\gamma)}{2}$, the term
$Z^{2n-l(\gamma)}_{2n-l-j}(...)$ is supported away from
$U^\gamma$. As the support of $Z^{2n-l(\gamma)}_{2n-l-j}(...)$ is
compact, the function $dist(\gamma^{-1}(x), x)$ achieves its
absolute minimum on $\text{supp}(Z^{2n-l(\gamma)}_{2n-l-j}(...))$,
which is strictly positive as $Z^{2n-l(\gamma)}_{2n-l-j}(...)$ is
supported away from $U^\gamma$. We assume that this minimum to be
$\alpha_0$. We prove that when $j>n-l(\gamma)/2$, the term
$\Psi^\gamma_j(\gamma^{-1}Z^{2n-l(\gamma)}_{2n-l(\gamma)-j}(\cdots)D_0,
\cdots, D_j)=0$.

According to its definition,
$\Psi^\gamma_j(\gamma^{-1}Z^{2n-l(\gamma)}_{2n-l(\gamma)-j}(\cdots)D_0,
\cdots, D_j)$ is computed by
\[
(-1)^{\frac{j(j+1)}{2}}\big[\int_{t\Delta_j}\Str(\gamma
Z^{2n-l(\gamma)}_{2n-l(\gamma)-j}(\cdots)D_0k_{s_0}^N[\bar{\partial}^*,
D_1]k_{s_1}^{N}\cdots [\bar{\partial}^*, D_k]k_{s_j}^N)ds_1\cdots
ds_j\big]_-.
\]
The term $\Str(\gamma
Z^{2n-l(\gamma)}_{2n-l(\gamma)-j}(\cdots)D_0k_{s_0}^N[\bar{\partial}^*,
D_1]k_{s_1}^{N}\cdots [\bar{\partial}^*, D_k]k_{s_j}^N)$ can be
computed by the following integral
\begin{equation}
\label{eq:integral} \int_{U^{\times
(j+1)}}\text{str}(\gamma^*_{\gamma^{-1}(x_0)}Z^{2n-l(\gamma)}_{2n-l(\gamma)-j}(\cdots)
(x_0))k^N_{s_0}(\gamma^{-1}(x_0),
x_1)\prod_{i=1}^j[\bar{\partial}^*, D_i]_{x_i}k_{s_i}^N(x_i,
x_{i+1})dx_0\cdots dx_j,
\end{equation}
where $x_{j+1}$ is $x_0$. We notice that the truncated heat kernel
$k_{s_i}^N(x_i, x_{i+1})$ can be chosen to be supported within the
open set $dist(x_i, x_{i+1})<\alpha_0/j+1$. Therefore, the product
\[
k^N_{s_0}(x_0, x_1)\prod_{i=1}^j[\bar{\partial}^*,
D_i]_{x_i}k_{s_i}^N(x_i, x_{i+1})\]
is supported in the following
open set
\[
dist(\gamma^{-1}(x_0),x_1)<\alpha_0/j+1,\ dist(x_i,
x_{i+1})<\alpha_0/j+1,\ i=1, \cdots, j-1,\ dist(x_j,
x_0)<\alpha_0/j+1,
\]
which is a subset of the following open set
\[
dist(\gamma^{-1}(x_0),x_0)<dist(\gamma^{-1}(x_0),x_1)+
dist(x_1,x_2)+\cdots+dist(x_j,x_0)<\alpha_0.
\]
But we know that $Z^{2n-l(\gamma)}_{2n-l(\gamma)-j}(\cdots)$ is
supported within the closed subset $dist(\gamma^{-1}(x),x)\geq
\alpha_0$. Therefore, we conclude that the integral
(\ref{eq:integral}) vanishes. And there is no contribution in the
sum (\ref{eq:pairing}) for $j>n-l(\gamma)/2$.

We conclude that the only nonzero contribution in the sum
(\ref{eq:pairing}) is from $j=n-l(\gamma)/2,\ k=2n-l(\gamma)$.
Next, we compute the contribution of these terms.\\

We have the following expression for
$Z_{n-l(\gamma)/2}^{2n-l(\gamma)}$, for
$\partial_i=\partial_{y^i}$,
\[
\begin{split}
&Z^{2n-l(\gamma)}_{n-l(\gamma)}(\partial_{1},\cdots,
\partial_{n-l(\gamma)/2}; \varphi_0, \cdots,
\varphi_{2n-l(\gamma)})\\
=&\varphi_0\frac{\partial \varphi_1}{\partial y^1}\cdots
\frac{\partial \varphi_{n-l(\gamma)/2}}{\partial
y^{n-l(\gamma)/2}}\bar{\partial}\varphi_{n-l(\gamma)/2+1}\cdots
\bar{\partial}\varphi_{2n-l(\gamma)}+\text{shuffles}.
\end{split}
\]

By the above arguments, we know that components of
$\bar{\partial}\varphi_{n-l(\gamma)/2+1}\cdots
\bar{\partial}\varphi_{2n-l(\gamma)}$ containing differentials
along normal directions to $U^\gamma$ have no contribution to
$\Psi^\gamma_{n-l(\gamma)/2}(\gamma\cdots)$ because their supports
are inside the closed set $dist(\gamma^{-1}(x),x)\geq \alpha_0$.
Therefore
\[
\begin{split}
&\sigma^\gamma_{2n-l(\gamma)}(c^\gamma_E(U); \varphi_0\otimes
\cdots \otimes \varphi_{2n-l(\gamma)})\\
=&(-1)^{(n-l(\gamma))(n-l(\gamma)+1)/2}\sum_{\pi \in
S_{n-l(\gamma)/2}}\Psi_{n-l(\gamma)/2}^\gamma(\gamma B_{2n-l(\gamma)},
y^{\pi(1)}, \cdots, y^{\pi(n-l(\gamma)/2)}),
\end{split}
\]
where $B_{2n-l(\gamma)}$ is a multiplication operator
\[
\begin{split}
B_{2n-l(\gamma)}=&\sum_{\pi\in
S_{2n-l(\gamma)}}\text{sign}(\pi)\varphi_0\frac{\partial
\varphi_{\pi(1)}}{\partial y^1}\cdots \frac{\partial
\varphi_{\pi(n-l(\gamma)/2)}}{\partial
y^{n-l(\gamma)/2}}\\
&\cdot \frac{\partial \varphi_{\pi(n-l(\gamma)/2+1)}}{\partial
\bar{y}^{n-l(\gamma)/2+1}}\cdots \frac{\partial
\varphi_{\pi(2n-l(\gamma))}}{\partial
\bar{y}^{2n-l(\gamma)}}d\bar{y}^1\wedge \cdots \wedge
d\bar{y}^{n-l(\gamma)/2}.
\end{split}
\]

We define $B=\varphi_0\frac{\partial \varphi_1}{\partial
y^1}\cdots \frac{\partial \varphi_{n-l(\gamma)/2}}{\partial
y^{n-l(\gamma)/2}}\bar{\partial}\varphi_{n-l(\gamma)/2+1}\cdots
\bar{\partial}\varphi_{2n-l(\gamma)}$. And as all $\varphi_i$ and
$y^i$ are $\gamma$ invariant, $B$ is $\gamma$ invariant. We look
at the expression
\[
\begin{split}
&\Psi_{n-l(\gamma)/2}(B, y^1, \cdots
y^{n-l(\gamma)/2})\\
=&\text{sign}(\gamma)\int_{t\Delta_{n-l(\gamma)/2}}Str(\gamma Bk_{s_0}^N[\bar{\partial}^*,
y^1]k_{s^1}^N\cdots [\bar{\partial}^*,
y^{n-l(\gamma)/2}]k_{s_{n-l(\gamma)/2}}^N)ds_1\cdots
ds_{n-l(\gamma)/2},
\end{split}
\]
where $\text{sign}(\gamma)$ is equal to
$(-1)^{(n-l(\gamma)/2)(n-l(\gamma)/2+1)/2}$. As $B,
[\bar{\partial}^*, y^i]$ are all differential operator of order 0,
it is sufficient to compute the leading term as $t\to 0$.

When $t\to 0$, we are reduced to a neighborhood of the origin of
$\complex^n$ with the standard metric. And we have the following
formulas for the operators in coordinates
\[
\partial=\sum d\bar{y}^i\frac{\partial}{\partial \bar{y}^i}, \ \ \
\bar{\partial}^*=-\sum\frac{\partial}{\partial
y^i}\iota_{\frac{\partial}{\partial y^i}}, \ \ \
\Delta_{\bar{\partial}}=-\sum_{j=1}^n\frac{\partial ^2}{\partial
y^j\partial \bar{y}^j},
\]
where $\iota$ is the subtraction. And the heat kernel is
\[
k_t(y,y')=\frac{1}{(\pi t)^{n}}e^{-\frac{|y-y'|^2}{t}}.
\]
Furthermore, we notice $[\bar{\partial}^*,
y^i]=\iota_{\partial/\partial \bar{y}^i}$, which commutes with the
heat kernel. And the expression of $\Psi_{n-l(\gamma)/2}(B, y^1,
\cdots, y^{n-l(\gamma)/2})$ is simplified to
\[
\begin{split}
=&\text{sign}(\gamma)\int_{t\Delta_{n-l(\gamma)/2}}Str(\gamma\overline{B}k_t^N)ds_1\cdots
ds_{n-l(\gamma)/2}\\
=&\text{sign}(\gamma)Str(\gamma\overline{B}k_t^N)\int_{t\Delta_{n-l(\gamma)/2}}ds_1\cdots
ds_{n-l(\gamma)/2}\\
=&\frac{(-1)^{n-l(\gamma)/2}t^{n-l(\gamma)/2}}{(n-l(\gamma)/2)!\pi^{n-l(\gamma)/2}}\text{sign}(\gamma)Str(\gamma\overline{B}k_t^N),
\end{split}
\]
where $\overline{B}$ is defined by
$B=\overline{B}d\bar{y}^1\wedge\cdots \wedge
d\bar{y}^{n-l(\gamma)/2}$, and the first numerical factor in the
last expression is the volume of $t\Delta_{n-l(\gamma)/2}$.

By \cite{BeGeVe}[Theorem 6.11], if $\text{supp}(\varphi_0)\cap
\cdots \cap \text{supp}(\varphi_{2n-l(\gamma)})$ does not contain
any $\gamma$ fixed point, then $Str(\gamma \overline{B}k_t^N)$
converges to 0 as $t\to 0$, and when $\text{supp}(\varphi_0)\cap
\cdots \cap \text{supp}(\varphi_{2n-l(\gamma)})$ contains $\gamma$
fixed points, then as $t\to 0$
\[
Str(\gamma \overline{B}k_t^N)\rightarrow
\frac{\tr_E(\gamma)}{t^{n-l(\gamma)/2}}\int_{M^\gamma}\frac{1}{\det(1-\gamma^{-1})}\overline{B}
d|y^1|\cdots d|y^{n-l(\gamma)/2}|.
\]

Finally, as we notice that different order of $y^i$ does not
change the limit, we have the conclusion
\[
\begin{split}
&\sigma^\gamma_{2n-l(\gamma)}(c^\gamma_E(U); \varphi_0\otimes
\cdots \otimes \varphi_{2n-l(\gamma)})\\
=&(-1)^{(n-l(\gamma)/2)(n-l(\gamma)/2+1)/2}\sum_{\pi \in
S_{n-l(\gamma)/2}}\Psi_{n-l(\gamma)/2}^\gamma(\gamma B_{2n-l(\gamma)},
y^{\pi(1)}, \cdots, y^{\pi(n-l(\gamma))})\\
=&(n-l(\gamma)/2)!\Psi_{n-l(\gamma)/2}(B, y^1, \cdots
y^{n-l(\gamma)/2})\\
\rightarrow&\frac{(-1)^{n-l(\gamma)/2}\tr_E(\gamma)}{(2\pi
i)^{n-l(\gamma)/2}}\int_{M^\gamma}\frac{1}{\det(1-\gamma^{-1})}
\varphi_0 d \varphi_1\cdots d\varphi_{2n-l(\gamma)}.
\end{split}
\]
\end{proof}

In the above Proposition \ref{prop:pairing} (4), we proved that
$j$ can not be strictly greater than $n-l(\gamma)$. The same
arguments also prove the following corollary.
\begin{corollary}\label{cor:supp}
If $\text{supp}(\varphi_0)\cap\cdots\cap
\text{supp}(\varphi_{k})\cap U^\gamma=\varnothing$, then for any
$D_0, \cdots, D_k\in A$
\[
\sigma^\gamma_{k}(D_0, \cdots, D_k; \varphi_0\otimes \cdots
\varphi_{k})=0.
\]
\end{corollary}
\subsection{Local expression of $\gamma$-Lefschetz number}
In this subsection, we will use the results developed in the
previous subsection to prove that the first $\gamma$-twisted trace
is equal to the third one.\\

We first observe that as $\bar{\partial}$ is $\gamma$-invariant,
the same argument as \cite{enfe}[Proposition 4.1] proves that
\[
\sum_{i=0}^{n}(-1)^i \tr_{\Omega^{(0,i)}(M;E)}(\gamma
De^{-t\Delta_{\bar{\partial}}})
\]
is independent of $t$ and is equal to $\Tr_{1,\gamma}$.\\

Let $(\tilde{U}_i)$ be an open cover of $M$ chosen as follows. We
start with the open over $(U_i)$ of $M^\gamma$ as is chosen in
Section 3.2. Fix $\epsilon_0>0$. Let $T^\gamma$ be a $2\epsilon_0$
tubular neighborhood of $M^\gamma$ in $M$ with the map $\pi:
T^\gamma\to M^\gamma$. Define $\tilde{U}_i=\pi^{-1}(U_i)$.
$(\tilde{U}_i)$ forms a cover of the tubular neighborhood
$T^\gamma$. Then we extend $(\tilde{U}_i)$ to a open cover of $M$
by requiring that the extra open sets will not intersect with the
$\frac{3}{2}\epsilon_0$ neighborhood of $M^\gamma$. We choose
$(\varphi_i)$ to be a partition of unity subordinate to the cover
$(\tilde{U}_i)$ such that $(\varphi_i)$ restricts to become a
partition of unity on $M^\gamma$ subordinate to the cover $(U_i)$.
Furthermore, by choosing a proper cut-off function, we can require
that $\varphi_i$ to be $\gamma$ invariant and constant along the
direction orthogonal to $M^\gamma$ within the $\epsilon_0$
neighborhood of $M^\gamma$ if $M^\gamma\cap
\text{supp}(\varphi_i)\ne \varnothing$. And we require
$(\tilde{U}_i)$ to have the following properties for $i_0< \cdots
<i_k$,
\begin{enumerate}
\item $U_{i_0}\cap \cdots\cap U_{i_k}$ is either 0 or contractible;
\item If $k>>0$, then $U_{i_0}\cap\cdots \cap U_{i_k}$
is empty.\\
\end{enumerate}

We have the following localization property about
$\Tr_{1,\gamma}(D)$ for a differential operator $D$.

\begin{proposition}\label{prop:cut-off}
Let $D$ be a differential operator on $M$, which is not
necessarily holomorphic. If the support of $D$ is away from
$M^\gamma$, then
\[
\sum_i(-1)^i\left[\tr_{\Omega^{(0,i)}(M,E)}(\gamma
De^{-t\Delta_{\bar{\partial}}})\right]_{-}=0.
\]
\end{proposition}
\begin{proof}
According to \cite{BeGeVe}[Proposition 2.46], the kernel
$p_t(x,y)$ of $De^{-t\Delta_{\bar{\partial}}}$ has the following
asymptotic expansion as $t\to 0$
\[
||p_t(x,y)-h_t(x,y)\sum_{i=-m}^Nt^i\Psi_i(x,y)||=O(t^{N-n-m}),
\]
where $h_t(x,y)=(4\pi t)^{-n}\exp(-dist(x,y)^2/4t)\Psi(d(x,y)^2)$,
and $\Psi$ is a cut-off function, and $m$ is the order of the
operator $D$.

We observe that as $D$ is supported away from $M^\gamma$,
$\Psi_i(x,y)$'s support is away from $M^\gamma \times M\subset
M\times M$ for any $i$.

Now the $\gamma$ trace of $De^{-t\Delta_{\bar{\partial}}}$ is
computed by
\begin{equation}\label{eq:trace}
\begin{split}
&\int_M\tr(\gamma_xp_t(\gamma^{-1}(x),x))dx\\
=&\int_M\tr(\gamma_x)h_t(\gamma^{-1}(x),x)\sum_{i=-m}^\infty
t^i\Psi_i(\gamma^{-1}(x),x) dx.
\end{split}
\end{equation}
As $\Psi_i(x,y)$'s support is away from $M^\gamma\times M$,
$\Psi_i(\gamma^{-1}(x),y)$'s support is also away from
$M^\gamma\times M$, and therefore $\Psi_i(\gamma^{-1}(x),x)$'s
support is away from $M^\gamma$.

As $M$ is compact, we know that for each $i$, the support of
$h_t(\gamma^{-1}(x), y)\Psi_i(\gamma^{-1}(x), y)$ is a compact
subset of $M\times M$, and therefore the support of
$h_t(\gamma^{-1}(x),x)\Psi_i(\gamma^{-1}(x),x)$ is a compact set
of $M$. Accordingly on the support of
$h_t(\gamma^{-1}(x),x)\Psi_i(\gamma^{-1}(x),x)$, the function
$dist(\gamma^{-1}(x),x)$ reaches its absolute minimum. As the
support of $\Psi_i(\gamma^{-1}(x),x)$ is away from $M^\gamma$, we
know that there is a positive number $\epsilon$ such that the
support of $h_t(\gamma^{-1}(x),x)\Psi_i(\gamma^{-1}(x),x)$,
$\min(dist(\gamma^{-1}(x),x))=\epsilon>0$.

On the other hand, we know that the support of the heat kernel
$h_t(x,y)$ can be chosen to be arbitrarily close to the diagonal
in $M\times M$. This forces the function $h_t(\gamma^{-1}(x),x)$
to be supported in the neighborhood
$dist(\gamma^{-1}(x),x)<\epsilon$. Considering the above arguments
which show that the support of $\Psi_i(\gamma^{-1}(x),x)$ is
outside the open set $dist(\gamma^{-1}(x),x)<\epsilon$ for any
$i$, we conclude that the function
$h_t(\gamma^{-1}(x),x)\Psi_i(\gamma^{-1}(x),x)$ has to vanish for
any $i$. Therefore, we conclude that $\Tr_{1,\gamma}(D)=0$.
\end{proof}

\begin{proposition}\label{prop:trace-1}
The two $\gamma$-twisted traces are same, $\Tr_{1, \gamma}=\Tr_{3,\gamma}$.
\end{proposition}
\begin{proof} We compute $\Tr_{1,\gamma}$ using the following formula
\[
\sum_{j=1}^n (-1)^j[\tr_{\Omega^{(0,j)}(M,E)}(\gamma
De^{-t\Delta_{\bar{\partial}}})]_-.
\]
By Proposition \ref{prop:cut-off}, we can use a cut-off function
to localize $D$ to be supported within the $\epsilon_0$
neighborhood of $M^\gamma$ without changing the value of the above
$\gamma$-twisted trace.

By the partition of unity $(\varphi_i)$, we can express the above
$\gamma$-twisted trace by
\[
\sum_i\sum_{j=1}^n
(-1)^j[\tr_{\Omega^{(0,j)}(M,E)}(\gamma\varphi_iDe^{-t\Delta_{\bar{\partial}}})]_-.
\]
By the assumption that $D$ is supported with the $\epsilon_0$
neighborhood of $M^\gamma$, all $\varphi_i$'s having non-trivial
contributions in the above sum are from pullbacks of a partition
of
unity of $M^\gamma$.\\

Using the definition of the pairing $\sigma_0^\gamma$, we have
that
\[
\sum_{j=1}^n
(-1)^j[\tr_{\Omega^{(0,j)}(M,E)}(\gamma\varphi_iDe^{-t\Delta_{\bar{\partial}}})]_-=\sigma_0^\gamma(D_i;\varphi_i),
\ \ \ D_i=D|_{\tilde{U}_i}\in \calD_E(\tilde{U}_i).
\]

As the twisted Hochschild homology of $\calD_E$ has 0 cohomology
in degree 0, we have that $D_i=b^\gamma D_i^{(1)}$. Then by
Proposition \ref{prop:pairing}, we have that
\[
\begin{split}
\Tr_{1, \gamma}
&=\sum_i \sigma^\gamma_0(b^\gamma D_i^{(1)};\varphi_i)\\
&\text{Use Proposition \ref{prop:pairing}, (ii)}\\
&=\sum_i\sigma^\gamma_1(D_i^{(1)}; 1,\varphi_i)\\
&\text{Use the partition of unity}\\
&=\sum_j\sum_i\sigma^\gamma_1(D_i^{(1)}; \varphi_j, \varphi_i)\\
&=\sum_{i\ne j}\sigma_1^\gamma(D^{(1)}_i; \varphi_j, \varphi_i)+\sum_j\sigma_1^\gamma(D^{(1)}_j; \varphi_j, \varphi_j)\\
&\text{Write $\varphi_j=1-\sum_{i\ne j}\varphi_i$}\\
&=\sum_{i\ne j}\sigma^\gamma_1(D^{(1)}_i-D^{(1)}_j;\varphi_j,\varphi_i)\\
&=\sum_{i\ne j}\sigma_1^\gamma(\check{\delta}D^{(1)}_{j,i}; \varphi_j, \varphi_i).
\end{split}
\]
According to Corollary \ref{cor:supp}, we know that if
$\text{supp}(\varphi_i)\cap \text{supp}(\varphi_j)$ is away from
$M^\gamma$, then the pairing
$\sigma_1^\gamma(\check{\delta}D^{(1)}_{j,i};\varphi_j,\varphi_i)=0$.
Therefore, by the assumption on supports of the partition of unity
on $M$, we conclude that in the above sum, all the nontrivial
contributions are from $\varphi_i,\varphi_j$ which are pullbacks
of a partition of unity on $M^\gamma$.\\

We can repeat this computation by induction and have the following equality,
\[
\begin{split}
\Tr_{1, \gamma}(D)&=\sum_{i_0<\cdots <i_{2n-l(\gamma)}}\sigma_{2n-l(\gamma)}^\gamma(s^\gamma_{i_0, \cdots, i_{2n-l(\gamma)}}; \varphi_{i_0, \cdots, i_{2n-l(\gamma)}})\\
+&\sum_{i_0<\cdots <i_{2n-l(\gamma)+1}}\sigma_{2n+1-l(\gamma)}^\gamma
(\check{\delta}D^{(2n+1-l(\gamma))}_{i_0, \cdots, i_{2n+1-l(\gamma)}};\varphi_{i_0, \cdots, i_{2n+1-l(\gamma)}}),
\end{split}
\]
where  $\varphi_{i_0, \cdots, i_{2n-l(\gamma)}}=\sum_{\pi \in
S_{2n-l(\gamma)+1}}\text{sign}(\pi)\varphi_{i_{\pi(0)}}\otimes
\cdots \otimes \varphi_{i_{\pi(2n-l(\gamma))}}$, and
$\varphi_{i_0}, \cdots, \varphi_{i_{2n-l(\gamma)}}$ are from
pullbacks of a partition of unity on $M^\gamma$.

For the second term in the above equation since there is no
further nontrivial Hochschild cycles, we can continue the
induction step and have  for $k\geq 2n+1-l(\gamma)$
\[
\sum_{i_0<\cdots <i_{2n-l(\gamma)+1}}\sigma_{2n+1-l(\gamma)}^\gamma
(\check{\delta}D^{(2n+1-l(\gamma))}_{i_0, \cdots, i_{2n+1-l(\gamma)}};\varphi_{i_0, \cdots, i_{2n+1-l(\gamma)}})=
\sum_{i_0<\cdots <i_k}\sigma_{k}^\gamma(\check{\delta}D^{(k)}_{i_0, \cdots, i_{k}};\varphi_{i_0, \cdots, i_{k}}).
\]
When $k$ is large enough, we know that
$\cap_{i=0}^k\text{supp}(\varphi_i)=\varnothing$ and by
Proposition \ref{prop:pairing} (i), these terms vanish.

Hence, we have that
\[
\Tr_{1, \gamma}(D)=\sum_{i_0<\cdots
<i_{2n-l(\gamma)}}\sigma_{2n-l(\gamma)}^\gamma(s^\gamma_{i_0,
\cdots, i_{2n-l(\gamma)}}; \varphi_{i_0, \cdots,
i_{2n-l(\gamma)}}).
\]
As all $\varphi_{i_k}$ in the above sum are from pullbacks of a
partition of unity on $M^\gamma$, according to our assumption at
the beginning of this subsection,  all these $\varphi_{i_k}$s are
$\gamma$ invariant and constant along the normal direction to
$M^\gamma$ within $\epsilon_0$ distance. We can apply Proposition
\ref{prop:pairing} (iv) to evaluate
$\sigma_{2n-l(\gamma)}^\gamma(\cdots)$, and have
\[
\Tr_{1,
\gamma}(D)=(2n+1-l(\gamma))!\frac{\tr_E(\gamma)(-1)^{n-l(\gamma)/2}}{(2\pi
i)^{n-l(\gamma)/2}}\sum_{i_0<\cdots <
i_{2n-l(\gamma)}}\lambda_{i_0, \cdots,
i_{2n-l(\gamma)}}(D)\int_{M^\gamma}\varphi_{i_0}d\varphi_{i_1}\cdots
d\varphi_{i_{2n-l(\gamma)}}.
\]
Since the restriction of $(\varphi_i)$ to $M^\gamma$ forms a partition of unity subordinate to $(U_i)$, we can evaluate the integral
\[
\int_{\Delta_{2n-l(\gamma)}}\varphi_0d\varphi_1\cdots d\varphi_{2n-l(\gamma)}=\frac{1}{(2n-l(\gamma)+1)!}.
\]
Hence we have that
\[
\begin{split}
\Tr_{1,\gamma}&=\frac{\tr_E(\gamma)(-1)^{n-l(\gamma)/2}}{(2\pi
i)^{n-l(\gamma)/2}}\sum_{i_0<\cdots
<i_{2n-l(\gamma)}}\lambda_{i_0, \cdots,
i_{2n-l(\gamma)}}(D)\text{sign}(i_0, \cdots,
i_{2n-l(\gamma)})\\
&=\frac{\tr_E(\gamma)(-1)^{n-l(\gamma)/2}}{(2\pi
i)^{n-l(\gamma)/2}}\Tr_{3,\gamma}.
\end{split}
\]
\end{proof}

In summary, we have proved the following formula for the
$\gamma$-twisted Lefschetz number.
\begin{theorem}\label{thm:global}
Let $M$ be a compact complex manifold, and $\gamma$ be an element
of a compact group acting on $M$ preserving the complex structure,
and $E$ be a $\gamma$-equivariant complex vector bundle on $M$,
and $D$ be a differential operator acting on $E$. Then
\[
\sum_i(-1)^i\tr(\gamma H^i(D))=\int_{M^\gamma} \frac{1}{(2\pi i)^{n-l(\gamma)/2}}\chi_{0,\gamma}(D),
\]
where $\chi_{0,\gamma}$ is as defined in Section 2.
\end{theorem}
\section{Orbifold Lefschetz number}
Let $G$ be a compact group acting on a complex manifold $M$ with
finite isotopy subgroups. Assume that the quotient $X=M/G$ to be
compact. In this section, we want to provide a $G$-equivariant
Lefschetz number formula for $G$-invariant differential operator
$D$ acting on a $G$-equivariant vector bundle $E$.  As $G$ acts on
the sheaf cohomology $H^j(M,E)$, we denote $H^j_G(M,E)$ to be the
subspace of $G$-fixed points in $H^j(M.E)$. Let $p$ be the
projection from $M\to X$. We  consider the pushforward vector
bundle $p_*(E)$ on $X$.  $H^j_G(M,E)$ can be identified with $H^j
(X, p_*(E))$. Because $D$ is $G$-invariant, $D$ acts on
$H^j_G(M,E)$. Furthermore, as $X$ is compact, we know that $H^j_G(M,E)$ is
finite dimensional. Therefore we can define the $G$-equivariant Lefschetz
number of $D$ to be
\[
\sum_i(-1)^i \tr(H^i_G(D)).
\]

We will use the geometric data on the orbifold to compute the
Lefschetz number.  Therefore, we recall some differential geometry
on a complex orbifold $X$. Given an orbifold $X$, we have a
naturally associated orbifold $\tilde{X}$, which is usually called
the inertia orbifold for $X$. Let us define this inertia orbifold
locally. In a sufficiently small open subset $U$ of $X$, we can
represent $X$ by a global quotient $V/\Gamma$, where $V$ is an
open subset of $\complex^n$ and $\Gamma$ is a finite group acting
on $V$ linearly. Accordingly, we introduce the following
stratified space $\tilde{U}$
\[
\coprod_{(\gamma)\subset \Gamma}V^\gamma/C(\gamma),
\]
where $(\gamma)$ stands for the conjugacy class of $\gamma$, and
$C(\gamma)$ is the centralizer group of $\gamma$ in $\Gamma$, and
$V^\gamma$ is the $\gamma$ fixed point subspace of $V$. The
stratified charts $\tilde{U}$ glue together to become a stratified
complex orbifold, which is usually denoted by $\widetilde{X}$. We
look at the $\chi_{0,\gamma}(D)$ as defined in Section 2.  The
pushforward of the collection $\chi_{0,\gamma}$ for all $\gamma$
from $V^\gamma$ to $V^\gamma/C(\gamma)$ forms a section of top
forms on the inertia orbifold $\widetilde{X}$. We remark that
because the quotient map $V^\gamma\to V^\gamma/C(\gamma)$ is a
proper locally embedding map, the pushforward map is well defined.

Let us first consider the situation that $G$ is a finite group.
\begin{theorem}\label{thm:finite}
Let $G$ be a finite group, and $\chi_{0,\gamma}(D)$ be the form defined in Section 2. We have
\[
\sum_i(-1)^i \tr(H^i_G(D))=\sum_{(\gamma)\subset
G}\int_{M^\gamma/C(\gamma)}\frac{1}{m_\gamma}\frac{1}{(2\pi
i)^{n-l(\gamma)/2}}\chi_{0,\gamma}(D),
\]
where $m_\gamma$ is the order of $\gamma$, a local constant number
on $M^\gamma$, and $(\gamma)$ runs over all conjugacy classes of
$G$.
\end{theorem}
\begin{proof}
We consider the $G$ averaging operator $P_G=\frac{1}{|G|}\sum_{\gamma\in G}\gamma$ acting on $H^i(X, E)$, where $|G|$ is the size of $G$. It is straightforward to check that $P_G$ is a projection from $H^i(X,E)$ to $H^i_G(X,E)$. Therefore, we have
\[
\begin{split}
\tr(H^i_G(D))&=\tr(P_GH^i(D)P_G)\\
&=\tr(P_GH^i(D))\\
&=\sum_{\gamma\in G}\frac{1}{|G|}\tr(\gamma H^i(D)).
\end{split}
\]
And applying Theorem \ref{thm:global}, we obtain the following equality for the Lefschetz number
\[
\begin{split}
\sum_i(-1)^i\tr(H_G^i(D))&=\sum_{\gamma\in G}\frac{1}{|G|}\sum_{i}(-1)^i\tr(\gamma H^i(D))\\
&=\sum_{\gamma\in G}\frac{1}{|G|}\int_{M^\gamma}\frac{1}{(2\pi i)^{n-l(\gamma)/2}}\chi_{0,\gamma}(D)\\
&=\sum_{(\gamma)\subset G}\frac{1}{|G|}\sum_{\alpha\in
(\gamma)}\int_{M^\alpha}\frac{1}{(2\pi
i)^{n-l(\alpha)/2}}\chi_{0,\alpha}(D).
\end{split}
\]
We notice that for different $\alpha,\alpha'$ in the same
conjugacy class of $G$, $M^\alpha$ is diffeomorphic to
$M^{\alpha'}$, $l(\alpha)=l(\alpha')$ and
$\chi_{0,\alpha}(D)=\chi_{0,\alpha'}(D)$ as $D$ and $\chi_{0,
\cdot}$ are both $G$-invariant. Let $|(\gamma)|$ and $|G|$ be the
sizes of $(\gamma)$ and $G$.  We continue the above computation
\[
\begin{split}
\sum_i(-1)^i\tr(H_G^i(D))&=\sum_{(\gamma)\subset G}\int_{M^\gamma}
\frac{|(\gamma)|}{|G|}\frac{1}{(2\pi i)^{n-l(\gamma)/2}}\chi_{0,\gamma}(D)\\
&=\sum_{(\gamma)}\int_{M^\gamma}\frac{1}{C(\gamma)}\frac{1}{(2\pi i)^{n-l(\gamma)/2}}\chi_{0,\gamma}(D)\\
&=\sum_{(\gamma)\subset
G}\int_{M^\gamma/C(\gamma)}\frac{1}{m_\gamma}\frac{1}{(2\pi
i)^{n-l(\gamma)/2}}\chi_{0,\gamma}(D),
\end{split}
\]
where in the last line of the above computation, we have used the definition of integration over an orbifold.
\end{proof}

In the last part of this section we consider a compact complex
orbifold $X$. Let $\calO_X$ be the sheaf of holomorphic functions
on $X$. By the sheaf of geometric holomorphic differential
operator on $X$, we mean the module of $\calO_X$ generated by
sections of the sheaf of holomorphic vector fields on $X$. A
global geometric differential operator $D$ on $X$ is a section of
the sheaf of geometric holomorphic differential operators on $X$.
Let $E$ be an orbifold vector bundle on $X$. We can define
geometric differential operators on $E$ in the same fashion. As
$X$ is compact, the generalized Lefschetz formula is defined as
same as on a manifold to be
\[
\sum_{i}(-1)^i\tr(H^i(D)).
\]
On the other hand, locally an orbifold is like a quotient of $\complex^n$ by a finite linear group action.  We can apply the construction of $\chi_0$ as in Theorem \ref{thm:finite} to $D$ and $\chi_0(D)$ is a differential form on the inertia $\widetilde{X}$.

\begin{theorem}\label{thm:orbifold}
Let $D$ be a global geometric differential operator acting on a
vector bundle $E$ of a compact complex orbifold $X$. Then
\[
\sum_{i}(-1)^i\tr(H^i(D))=\int_{\widetilde{X}}\frac{1}{m_\calO}\frac{1}{(2\pi i)^{n-l(\calO)/2}}\chi_{0, \calO}(D),
\]
where $m_{\calO}$ is a local constant on $\widetilde{X}$ telling the size of isotopy, $l(\calO)$ is a local constant telling the codimension of $\calO$ inside $X$, and $\chi_{0, \calO}$ is a top degree differential form on $\widetilde{X}$.
\end{theorem}
\begin{proof}
We start with the following local result from \cite{AFLS:HIAWAGF}.
For any $x\in X$, there is a small enough neighborhood $U$ which
is holomorphic to the quotient of a complex open set $V$ by a
finite group $\Gamma$ linear action. And the Hochschild
(co)homology of $\calD_E(U)$ is spanned by the conjugacy classes
of $\Gamma$ with degree equal to the codimension of the fixed
point subspace.  A crucial observation to the proof of this result
is that the $\Gamma$-invariant subalgebra of the Weyl algebra is
Morita equivalent to the crossed product algebra of the Weyl
algebra with $\Gamma$. The generators of (co)homology are the sums
of  $c^\gamma_{2n-l(\gamma)}$ and $\tau^\gamma_{2n-l(\gamma)}$ in
the same conjugacy classes.

The proof of this theorem is essentially a copy of Section 3. We
can choose a nice cover of $X$ and the corresponding
triangulation.  The new ingredient is that we will have more than
one contributions in the arguments of Proposition
\ref{prop:trace-2}-\ref{prop:trace-1}, and instead we will have
one contribution for each conjugacy class. However, since we have
1-1 correspondence for each $\gamma$ by Proposition
\ref{prop:trace-2}-\ref{prop:trace-1}, by taking the sum, we
obtain the statement of this theorem.\\
\end{proof}

\begin{remark}
\label{rmk:diff} We remark that the class of geometric
differential operators considered in Theorem \ref{thm:orbifold} is
quite restrictive. It excludes many interesting operators which
should be viewed as differential operators in an algebraic way,
where were called algebraic differential operators in
Introduction.

For example, we consider the simplest complex orbifold
$\complex/\Z_2$, where $\Z_2$ acts on $\complex$ by $z\mapsto -z$.
The differential operators considered in Theorem
\ref{thm:orbifold} are all of the form $\sum_i f_i\partial_z^i$,
where $f_i$ is an even polynomial when $i$ is even, and an odd
polynomial otherwise.  These are $\Z_2$ invariant differential
operators on $\complex$.

On the other hand, the algebra of polynomials on $\complex$
invariant under $\Z_2$ are even polynomials. The operator
$D=\frac{1}{z}\partial_z$ acts on the space of even polynomials
linearly satisfying the Leibniz rule
\[
\frac{1}{z}\partial_z(z^{2m})=2mz^{2m-2}, \ \ \ \forall m\geq0.
\]
But it is obvious that this operator cannot be written as
$\Z_2$-invariant differential operator on $\complex$. This
operator is a $\Z_2$-invariant meromorphic differential operator
on $\complex$ which descends to a ``differential operator" on
$\complex/\Z_2$.  This is probably related to connections with
logarithmic singularities \cite{langer}. A Lefschetz formula for
this type of operator is a future research project.
\end{remark}

\bibliographystyle{alpha}

\vskip 2mm
\address{
\noindent{Giovanni Felder}, {\tt felder@math.ethz.ch}\newline
    \noindent{\rm Department of Mathematics, ETH Zurich, 8092 Zurich,
    Switzerland}\\

\noindent{Xiang Tang}, {\tt xtang@math.wustl.edu} \newline
   \noindent {\rm  Department of Mathematics, Washington University, St. Louis,
           USA }
}

\end{document}